\documentclass[12pt]{article}
\usepackage{amssymb,amsfonts}
\usepackage{amsfonts}
\usepackage{latexsym}
\usepackage[dvips]{graphics}
\usepackage{mathdots,color}
\usepackage{bm}
\usepackage{amsmath}
\usepackage{mathrsfs}
\usepackage{epsf, epsfig, color, epstopdf}

\newtheorem{lem}{Lemma}
\newtheorem{theo}{Theorem}
\newtheorem{cor}{Corollary}
\newtheorem{prob}{Problem}
\newtheorem{cla}{Claim}
\newtheorem{prop}{Property}

\newcommand{\proof}{{\noindent {\em Proof}.\quad}\setcounter{countclaim}{0}
\setcounter{countcase}{0}}
\newcommand{\proofend}{{\hfill$\Box$}}

\newcounter{countfig}

\newcounter{countclaim}
\def\inclaim{\addtocounter{countclaim}{1}
{\noindent {\bf Claim \thecountclaim}: }}

\newcounter{countcase}

\newcommand{\beeq}{\begin{equation}}
\newcommand{\eneq}{\end{equation}}

\newcommand{\beeqn}{\begin{eqnarray*}}
\newcommand{\eneqn}{\end{eqnarray*}}

\newcommand\red[1] {{\color{red} #1}}
\newcommand\green[1] {{\color{green} #1}}
\newcommand\blue[1] {{\color{blue} #1}}

\def \usecolour
{\newcommand {\rered}{\red}
\newcommand {\reblue} {\blue}
\newcommand {\regreen} {\green}
}

\usecolour   
\newcommand {\relabel}[1] {\label{#1} \red{[*: #1]}}

\def\relabel {\label}   

\begin{document}

\newcommand{\resection}[1]
{\section{#1}\setcounter{equation}{0}}

\renewcommand{\theequation}{\thesection.\arabic{equation}}

\baselineskip 0.6 cm

\title {A new infinite family of 4-regular
crossing-critical graphs
}

\author
{Zongpeng Ding\thanks{Email: dzppxl@163.com.}\\
\small Department of Mathematics\\
\small Hunan First Normal University,
Changsha 410205, China\\
\\
Yuanqiu Huang\thanks{Corresponding author. Email: hyqq@hunnu.edu.cn.}\\
\small  School of Mathematics and Statistics\\
\small Hunan Normal University, Changsha 410081, China\\
\\
Fengming Dong\thanks{Email: donggraph@163.com
and fengming.dong@nie.edu.sg.}\\
\small National Institute of Education,  Nanyang Technological
University,  Singapore
}

\date{}

\maketitle

\begin{abstract}
A graph $G$ is said to be crossing-critical if $cr(G-e)< cr(G)$ for every edge $e$ of $G$, where $cr(G)$
is the crossing number of $G$.
Richter and Thomassen [Journal of Combinatorial Theory, Series B 58 (1993), 217-224] constructed an infinite family of 4-regular crossing-critical
graphs with crossing number $3$.
In this article,  we present a new infinite family of 4-regular crossing-critical
graphs.

\end{abstract}




\resection{Introduction}
All graphs considered here are undirected, finite and simple.
For any terminology and definition without explanation, we refer to \cite{BA}.
For a graph $G$,
let $V(G)$ and $E(G)$ denote its
vertex set  and edge set
respectively.
For any edge set $E_0$ of $G$,
let $G\setminus E_0$ denote the spanning
subgraph of $G$ obtained by deleting
all edges in $E_0$.
In particular, if $E_0=e$,
then $G\setminus E_0$ is written as $G\setminus e$.

For a drawing $D$ of $G$ in the plane,
let $cr_{D}(G)$ denote the number of crossings of edges in $D$.
The {\it crossing number} of  $G$, denoted by $cr(G)$, is the
minimum value of $cr_{D}(G)$
over all drawings $D$ of $G$ in the plane.
A graph $G$ is said to be
{\it crossing-critical} if $cr(G\setminus e)< cr(G)$ for every edge $e$ of $G$.

The crossing number of a graph is a parameter that measures how far a graph is
from a planar graph, and is a classic topological invariant of the graph. Its
theory has been widely applied to the repaint and identification of sketch, layout
problem in the VLSI in large scale circuit, and the graphical representation of
DNA in biology engineering and so on
(see \cite{CM,LA,NS}).

Computing the crossing number of a graph is an NP-hard problem \cite{GM}.
For example,
$cr(K_{7,11})$ and $cr(K_{13})$ are till unknown.
There are still many graph classes whose crossing numbers have not be confirmed
(see \cite{J} and \cite{SM}).
For the initial work about the crossing numbers of the complete bipartite graphs $K_{m,n}$ and complete graphs $K_{n}$, one can refer to \cite{BL}.

For $K_{m,n}$, Zarankiewicz
proposed the followingc conjecture \cite{Z}:
$$
cr(K_{m,n})
=\left\lfloor\frac{m}{2}
\right\rfloor
\left\lfloor\frac{m-1}{2}\right\rfloor\left\lfloor\frac{n}{2}\right\rfloor\left\lfloor\frac{n-1}{2}\right\rfloor.
$$
To date, the conjecture has been verified for $\min\{m, n\}\leq 6$
due to Kleitman~\cite{K}
and for $m=7$ and $n \leq 10$
due to Woodall~\cite{W}.

For $K_{n}$,
Guy  \cite{GR} conjectured that
$$cr(K_n)=\frac{1}{4}
\left\lfloor\frac{n}{2}\right\rfloor
\left \lfloor\frac{n-1}{2}\right\rfloor
\left\lfloor\frac{n-2}{2}\right\rfloor
\left\lfloor\frac{n-3}{2}\right\rfloor.
$$
Up to now, the conjecture has been
verified for $n \leq 12$
by McQuillan and Richter \cite{M}.
McQuillan, Pan  and Richter
\cite{MD} has also proved that
$cr(K_{13})\ne 217 $.

It is even more challenging
to determine the crossing numbers of cross-critical graphs.
For the  extensive literature on infinite families of crossing-critical graphs,
we may refer to
\cite{B,H,RB,S}.
Richter and Thomassen \cite{RB} constructed an infinite family of 4-regular crossing-critical
graphs with crossing number $3$, and posed the following problem.

\begin{prob}\label{prob1}
Are there infinitely many 5-regular crossing-critical
graphs?
\end{prob}

In this article,
we will present  a new infinite family of 4-regular crossing-critical
graphs $G_{n}$ of order $2n$ and
crossing number $cr(G_n)=n$ for $n\ge 4$.

\begin{figure}[h!]
\centering
\includegraphics[width=9cm]
{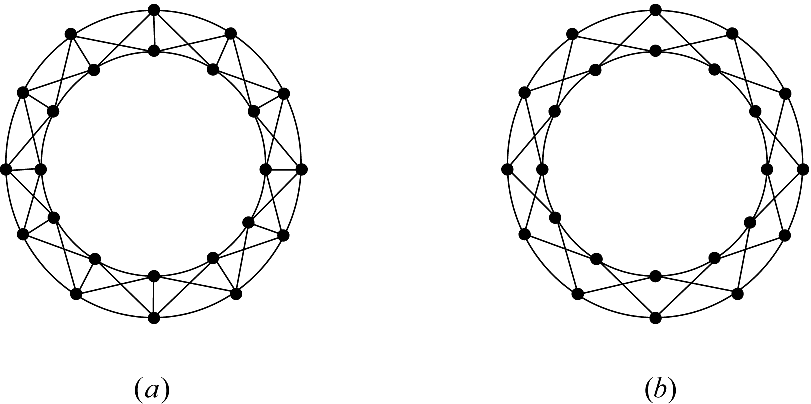}
\caption{(a) Musical graph and (b)  $G_{12}$.}\label{fig1}
\end{figure}

The graph shown in Figure  \ref{fig1} (a)
was called the {\it musical graph}
by Donald E. Knuth  (\cite{knuth},
Problem 133 in p. 44).
While many properties of this graph can be easily analyzed,
the question ``Can it be drawn with
fewer than 12 crossings?"
was finally settled by
Petra Mutzel~\cite{PM} who used a
program to show that
the crossing number of
this graph is indeed equal to $12$.

The musical graph, showed in
Figure~\ref{fig1} (a),
is actually
the graph $C_{12}\circ K_2$,
where
$G_1\circ G_2$,
called the
{\it lexicographic product}
of $G_{1}$ and $G_{2}$,
is the graph such that
\begin{itemize}
	\item its vertex set is
the cartesian product
$V(G_{1})\times V(G_{2})$; and
\item
any two vertices $(u, u^{\prime})$ and $(v, v^{\prime})$
are adjacent in $G_1\circ G_2$
if and only if either $u$ is adjacent to  $v$ in $G_{1}$ or $u = v$ and $u^{\prime}$ is
adjacent to $v^{\prime}$ in $G_{2}$.
\end{itemize}

Let $M_n$ denote the graph
$C_{n}\circ K_2$.
Then $M_{12}$ is the musical graph
and  Mutzel \cite{PM} showed that
$cr(M_{12})=12$.
Ouyang et al. ~\cite{OZ}
extended the conclusion
to $cr(M_{n})=n$ for all $n\ge 3$.

For $n\ge 3$,
let $G_{n}$ denote the lexicographic product $C_n \circ \overline{K_2}$.
For example, $G_{12}$ is
the graph shown in Figure
\ref{fig1} (b).
Obviously, $G_{n}$ is a spanning
subgraph of $M_{n}$
for all $n\ge 3$.
Note that $G_{3}$ is planar
while $G_n$ is not for all $n\ge 4$.
The main purpose in this article is to show that for
$n\ge 4$, $G_n$ has crossing number $n$ and $G_n$ is crossing-critical.
\begin{theo}\label{theo1}
	For any $n\ge 4$,  $cr(G_n)=n$ and $G_n$ is crossing-critical.
\end{theo}

Obviously, the result of Theorem~\ref{theo1}
implies that $cr(M_n)=n$
for all $n\ge 4$.


\begin{cor}[\cite{PM, OZ}]
	\label{cor1}
	For any integer $n\ge 3$,
	$cr(M_n)=n$.
\end{cor}

\resection{Proof of Theorem~\ref{theo1}}

Note that $G_n$ is actually the graph
with vertex set $\{x_i,y_i:1\le i\le n\}$
and edge set
$$
\{x_ix_{i+1},x_iy_{i+1},y_ix_{i+1},y_iy_{i+1}:i=1,2,\ldots, n\},
$$
where $x_{n+1}=x_1$ and $y_{n+1}=y_1$.

\begin{figure}[htbp]
\centering
\includegraphics[width=5 cm]
{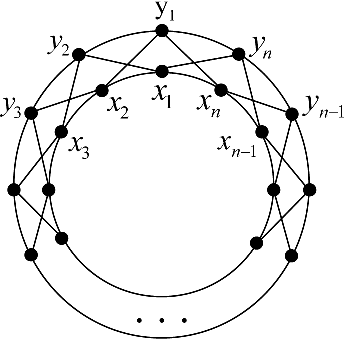}
\caption{The graph $G_{n}\cong C_n \circ \overline{K_2}$.}\label{fig2}
\end{figure}

Figure \ref{fig2},
$G_{n}$ is drawn in the plane
with $n$ crossings.
This indicates that $cr(G_{n})\leq n$.

\noindent {\it Proof of Theorem~\ref{theo1}}:
We first show that
 $cr(G_{n})=n$ for all $n\ge 4$.
Note that $G_{4}\cong K_{4,4}$,
implying that
$cr(G_{4})=cr(K_{4,4})=4$
due to Kleitman \cite{K}.
Suppose now that $n\ge 5$ and
	$cr(G_{k})=k$ holds for $4\leqslant k\leqslant n-1$.


Figure \ref{fig2} shows that $cr(G_{n})\leq n$.
Thus, it suffices to show that
$cr_{\phi}(G_{n})\geq n$ in any drawing $\phi$.

For $1\le i\le n$, let $L^i$ denote the
subgraph of $G_n$ induced by
$\{x_{i},x_{i+1},y_{i},y_{i+1}\}$,
where $x_{n+1}=x_1$
and $y_{n+1}=y_1$.
Clearly, $E(L^i)=\{x_ix_{i+1},x_iy_{i+1},
y_{i}x_{i+1}, y_{i}y_{i+1}\}$.
For $1\leq i<j\leq n$, $E(L^{i}) \cap E(L^{j})=\emptyset$,  and $V(L^{i}) \cap V(L^{j})\ne \emptyset$ if and only if
either $j=i+1$ or $j=n$ and $i=1$.

In the following,
let $\gamma_{\phi}(L^{i})$ be
number of edges in $E(L^i)$
which are crossed in $\phi$.
We are now going to accomplish the proof by showing the
following claims.

\begin{cla}\label{cl1}
If $\gamma_{\phi}(L^{i})=1$ holds
for some $i$ with $1\leq i\leq n$,
then $cr_{\phi}(G_n)\ge n$.
\end{cla}

Suppose that $\gamma_{\phi}(L^{1})=1$.
Clearly, edges in $L^{1}$  do not generate any crossings in $\phi$.
Otherwise, $\gamma_{\phi}(L^{1})\geq2$.
This indicates that $cr_{\phi}(L^{1})=0$.
Without loss of generality,
let $x_{1}y_{2}$ be the only crossing edge in $L^{1}$.
We remove the edges $x_{1}y_{2}$ and $x_{2}y_{1}$ from $\phi$
and obtain the restricted drawing $\phi_{1}$ of the subgraph $G_{n}\setminus\{x_{1}y_{2},x_{2}y_{1}\}$ induced by $\phi$.
Note that $x_{1}y_{2}$ is a crossing edge, then at least one crossing was reduced by removing these two edges.
Thus,
$$ cr_{\phi}(G_{n})\geq cr_{\phi_{1}}
(G_{n}\setminus \{x_1y_2,x_2y_1\})
+1.
$$

In $\phi_{1}$, since $x_{1}x_{2}$ and $y_1y_2$ are not crossing edges,
by contracting both $x_{1}x_{2}$ and $y_{1}y_{2}$,
we get a restricted drawing $\phi_{2}$ of $G_{n-1}$ induced by $\phi$.
Obviously, the inequality $cr_{\phi_1}(G_{n}\setminus \{x_1y_2, y_1x_2\})
\ge cr_{\phi_2}(G_{n-1})$ holds.
Thus,
\begin{eqnarray*}
cr_{\phi}(G_{n})&\geq&cr_{\phi_{1}}(G_{n}\setminus\{x_{1}y_{2},x_{2}y_{1}\})+1\\
&=& cr_{\phi_{2}}(G_{n-1})+1\\
&\geq& n-1+1=n.
\end{eqnarray*}
Thus, Claim~\ref{cl1} holds.

	For any $U\subseteq \{1,2,\ldots,n\}$,
	let $\overline U=\{1,2,\ldots,n\}\setminus U$ and let $E_U$ denote the set
	$\bigcup_{i\in U}E(L^i)$.
	For any edge sets $E_1$ and $E_2$
	of $G_n$,
	let $cr_{\phi}(E_1,E_2)$
	denote the number of crossings
	in $\phi$ between edges in $E_1$
	and edges in $E_2$.

\begin{cla}\label{cl2}
For any $U\subseteq \{1,2,\ldots,n\}$, if $\gamma_{\phi}(L^i)\ge 2$ for each $i\in U$,  then
$$
cr_{\phi}(E_U, E_U)
+cr_{\phi}(E_U, E_{\overline U})
\ge |U|.
$$
	\end{cla}

	Since $\gamma_{\phi}(L^i)\ge 2$
	for each $i\in U$,
	in the drawing $\phi$, there are at least
	$2|U|$ crossing edges in the edge set
	$E_U$.
	Each crossing point involves only two crossing edges.
	Thus, Claim~\ref{cl2} holds.

For any integer $j$ with $0\le j\le n+1$,
let $\mu(0)=n$, $\mu(n+1)=1$
and $\mu(j)=j$ otherwise.

\begin{cla}\label{cl3}
	For $1\le i\le n$,
if $\gamma_{\phi}(L^{i})=0$, then $cr_{\phi}(E(L^{\mu(i-1)}),E(L^{\mu(i+1)}))\geq4$.
\end{cla}

\begin{figure}[htbp]
\centering
\includegraphics[width=0.5\textwidth]{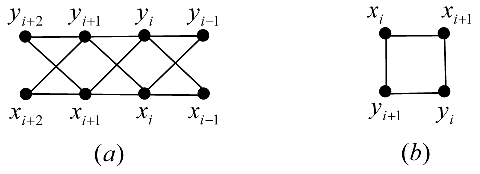}
\caption{(a) $L^{i-1}\cup L^{i}\cup L^{i+1}$ and  (b) $\phi(L^{i})$.}\label{figb}
\end{figure}

Without loss of generality, assume that
$\gamma_{\phi}(L^{i})=0$,
where $2\le i\le n-1$.
We are now going to show that
$cr_{\phi}(E(L^{i-1}),E(L^{i+1}))\geq4$.

By the assumption, there are no crossing edges in $L^i$, as shown in
Figure \ref{figb} (b).
Obviously, the subgraph, denoted by $G'$,
obtained from $G_n$ by  removing all four vertices in  $\{x_{i},x_{i+1},y_{i},y_{i+1}\}$  is connected.
Thus, the vertices of $G'$ are placed in the same region of  $\phi(L^{i})$; otherwise, some edges in $L^{i}$
are crossed.
We may assume that all vertices in $G'$ are within
the finite region of $\phi(L^{i})$,
implying that
$cr_{\phi}(E(P),E(P'))\ge 1$ holds
for any pair of vertex-disjoint paths
$P$ and $P'$ of $G_n$, where
$P$ joins $x_i$ to $y_i$
and $P'$ joins $x_{i+1}$ to $y_{i+1}$.

There are two edge-disjoint paths $P_1$ and $P_2$ in $L^{i-1}$  joining $x_i$ to $y_i$, where
$P_{1}=x_{i}x_{i-1}y_{i}$ and $ P_{2}=x_{i}y_{i-1}y_{i}$,
and two edge-disjoint paths $P_3$ and $P_4$  in $L^{i+1}$ joining $x_{i+1}$ to $y_{i+1}$, where
 $P_{3}=x_{i+1}x_{i+2}y_{i+1}$
and $P_{4}=x_{i+1}y_{i+2}y_{i+1}$.
Thus,
$$
cr_{\phi}(E(P_{k}),E(P_{l}))\geq1
$$
for any $k=1,2$ and $l=3,4$.
Therefore,
\begin{eqnarray*}
cr_{\phi}(E(L^{i-1}),E(L^{i+1}))&=&cr_{\phi}(E(P_{1}\cup P_{2}),E(P_{3}\cup P_{4}))\\
&=& cr_{\phi}(E(P_{1}),E(P_{3}))+cr_{\phi}(E(P_{1}),E(P_{4}))\\
& &+cr_{\phi}(E(P_{2} ),E(P_{3}))+cr_{\phi}(E(P_{2} ),E(P_{4}))\\
&\geq &4.
\end{eqnarray*}
Claim~\ref{cl3} holds.

Let $U_0=\{1\le i\le n: \gamma_{\phi}(L^{i})=0\}$ and
$$
U_1=\{1,2,\ldots,n\}
\setminus \{\mu(i),\mu(i-1),\mu(i+1): i\in  U_0\}.
$$
Obviously,
$U_1\subseteq \overline {U_0}$ and
$3|U_0|+|U_1|\ge n$.

\begin{cla}\label{cl4}
$cr_{\phi}(G_n)\ge n$.
\end{cla}

 By Claim~\ref{cl1},
 the result holds whenever
 $\gamma_{\phi}(L^i)=1$ for some $i$ with $1\le i\le n$.
 By Claim~\ref{cl2},
 the result also holds when
$\gamma_{\phi}(L^i)\ge 2$ for each $i$ with $1\le i\le n$.
Thus, we need only consider the case
that $U_0\ne \emptyset$
and $\gamma_{\phi}(L^i)\ge 2$
for each $i\in \overline {U_0}$.
Obviously, $\gamma_{\phi}(L^i)\ge 2$
for each $i\in U_1$.
It follows that
\begin{eqnarray*}
	cr_{\phi}(G_{n})&\ge &
	cr_{\phi}(E_{U_1},E_{U_1})
	+cr_{\phi}\left (E_{U_1},E_{\overline {U}_1}\right )
	+	\sum_{j\in U_0}
	cr_{\phi}\left (E(L^{\mu(j-1)}),  E(L^{\mu(j+1)})\right )
	\\
	&\ge & |U_1|+4|U_0|
	\\
	&\ge &n,
\end{eqnarray*}
where the second last inequality follows from Claims~\ref{cl2} and~\ref{cl3}.
Claim~\ref{cl4} holds
and thus
 $cr(G_{n})=n$ for all $n\ge 4$.

We now show that for any $n\ge 4$,
$G_n$ is crossing-critical.
Let $e$ be any edge in $G_n$.
Without loss of generality,
assume that $e$ is an edge in $L^1$.
If $e\in \{x_1y_2, x_2y_1\}$,
then $G_{n}\setminus e$ has a drawing $D$ obtained
from the one in Figure~\ref{fig2} by removing $e$.
Note that $cr_{D}(G_n\setminus e)=n-1$, implying that
$cr(G_n\setminus e)\le n-1$.

\begin{figure}[htbp]
\centering
\includegraphics[width=0.6\textwidth]{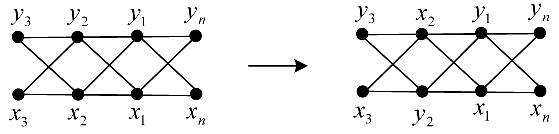}
\caption{Exchange
$x_2$ and $y_2$.}\label{fig9}
\end{figure}

If $e\in \{x_1x_2, y_1y_2\}$,
then we consider a new drawing of $G_n$ by simply exchanging
$x_2$ and $y_2$. See Figure~\ref{fig9}.
Similarly,
we can obtain a drawing of $G_{n}\setminus e$ in which the crossing number is less than $n$.
Thus, $G_n$ is crossing-critical for $n\ge 4$.
\proofend

\section{Concluding remarks}

 In \cite{RB},
 Richter and Thomassen
 asked if there are infinitely many 5-regular crossing-critical
graphs.
Let $H_n$ denote the Cartesian product $G_n\square  K_2$,
where
the {\it Cartesian product} $G\square  H$ of graphs $G$ and $H$ is the graph such that:
\begin{itemize}
	\item its vertex set is the Cartesian product $V(G) \times V(H)$; and
\item
two vertices $(u,v)$ and $(u' ,v' )$ are adjacent in $G \square H$ if and only if either
$u = u'$ and $v$ is adjacent to
$v'$ in $H$, or
$v = v'$ and
$u$ is adjacent to $u'$ in $G$.
\end{itemize}

\begin{figure}[h!]
	\centering
\includegraphics[width=8 cm]
{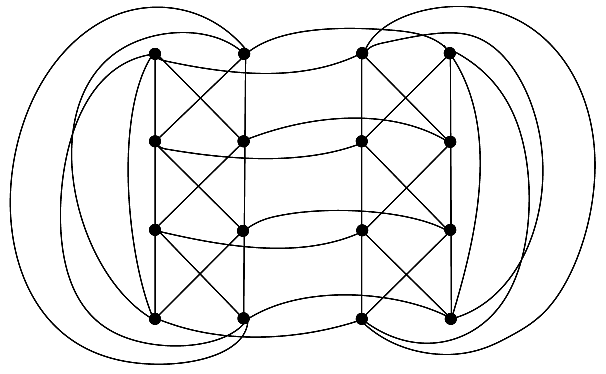}
	\caption{Graph $H_4=G_4\square K_2$.}\label{H4}
\end{figure}

Clearly, $H_n$ is 5-regular.
For example, $H_4$ is the graph shown in
Figure~\ref{H4}.
It is readily checked that
the graph $H_n$ has a  drawing
$D$ in the plane with $6n$ crossings.
We wonder if this drawing is optimal.

\begin{prob}
	For $n\ge 4$, is $H_n$ crossing-critical
	and  $cr(H_n)=6n$?
\end{prob}

\textbf{ Acknowledgments}

The first author would like to thank the hospitality of the  National Institute of Education, Nanyang Technological University in Singapore, where the work was done.

This work was supported by the National
Natural Science Foundation of China (No. 12371346, 12271157, 12371340),
the Scientific Research Fund of Hunan Provincial Education Department (No. 22A0637) and the Hunan Provincial Natural Science Foundation of
China (No. 2023JJ30178, 2022JJ30028).

\textbf{Data availability}

No data was used for the research described in the article.

\textbf{Conflict of interest}

We claim that there is no conflict of interest in our paper.


\begin{thebibliography}{99}\setlength{\itemsep}{-2mm}



\bibitem{B}D. Bokal, Infinite families of crossing-critical graphs with prescribed average degree and crossing number, {\it J. Graph Theory} {\bf 65}(2) (2010), 139-162.

\bibitem{BA} J. A. Bondy, U.S.R. Murty, Graph Theory, GTM 244, Springer, 2008.

\bibitem{BL} Lowell W. Beineke, Robin Wilson, The Early History of the Brick Factory
Problem, {\it Math. Int.}
{\bf 32}(2) (2010), 41-48.

\bibitem{CM}
M. Chimani, C. Gutwenger, Non-planar core reduction of graphs,
{\it Discrete Math.}
{\bf 309} (2009), 1838-1855.





\bibitem{GM}
M. R. Garey, D. S. Johnson, Crossing number is NP-complete,
{\it SIAM J. Algebraic Discrete Methods}
{\bf 4} (1983), 312-316.




\bibitem{GR}R. K. Guy, Crossing numbers of graphs, in Graph Theory and Applications, Y. Alavi, D. R. Lick, and A. T. White, eds., Lecture Notes Math. 303, Springer, New York, (1972) 111-124.

\bibitem{H} P. Hlineny, New infinite families of almost-planar crossing-critical graphs, {\it Electronic J. Combin.} {\bf 15(1)}   (2008), 1615-1615.

\bibitem{J}
P. K. Jha, S. Devisetty, Orthogonal drawings and crossing numbers of the Kronecker product of
two cycles,
{\it J. Parallel Distrib. Comput.}
{\bf 72} (2012), 195-204.

\bibitem{K} D. J. Kleitman, The crossing number of  $K_{5,n}$,
{\it J. Combin. Theory}
{\bf 9} (1970), 315-323.

\bibitem{knuth}
D. E. Knuth,  {\it The Art of Computer Programming: Combnatorial Algorithms, Part 1, Vol. 4A,}
Addison-Wesley, 2011.





\bibitem{LA}
A. Liebers, Planarizing graphs-A survey and annotated bibliography,
{\it J. Graph Algorithms Appl.}
{\bf 5} (2001), 1-74.

\bibitem{M} D. McQuillan, R. B. Richter, On the Crossing Number of without Computer Assistance,
{\it J. Graph Theory}
{\bf 82}(4) (2016), 387-432.


\bibitem {MD} D. McQuillan, S. Pan, R. B. Richter, On the crossing number of $K_{13}$,
{\it J. Combin. Theory Ser. B}
{\bf 115} (2015), 224-235.


\bibitem{PM}
Petra Mutzel, The crossing number of graphs: Theory and computation.
In {\it Efficient Algorithms}, pp. 305-317. Springer, Berlin, Heidelberg, 2009.

\bibitem{NS}
Sandeep N. Bhatt and F. Thomson Leighton, A framework for solving VLSI graph layout problems,
{\it J. Comput. System Sci.}
{\bf 28} (1984), 300-343.

\bibitem{OZ} Z. Ouyang, J. Wang and Y. Huang, The strong product of graphs and crossing numbers,
{\it Ars Combinatoria}
{\bf 137} (2018), pp. 141-147.




\bibitem{RB} R. B. Richter and C. Thomassen, Minimal graphs with crossing number at
least $k$,
{\it J. Combin. Theory Ser. B}
{\bf 58} (1993), 217-224.

\bibitem{S} G. Salazar, Infinite families of crossing-critical graphs with given average degree, {\it Discrete Math.} {\bf 271}(1-3) (2003), 343-350.

\bibitem{SM} M. Schaefer, Crossing numbers of graphs, CRC Press, Florida,
2017.

\bibitem{W} D. R. Woodall, Cyclic-order graphs and zarankiewicz's crossing number conjecture,
{\it J. Graph Theory}
{\bf 17} (1993), 657-671.


\bibitem{Z} K. Zarankiewicz, On a problem of P. Tur\'on concerning graphs, {\it Fund. Math} {\bf 41} (1955), 137-145.
\end{thebibliography}
\end{document}